\documentclass[11pt]{article}

\usepackage{amsfonts,latexsym,amsmath,amscd,geometry}
\usepackage{amssymb}
\usepackage{xcolor}
\usepackage[title]{appendix}

\usepackage[colorlinks=true, allcolors=brown]{hyperref}

\newtheorem{theorem}{Theorem}[section]
\newtheorem{lemma}[theorem]{Lemma}
\newtheorem{proposition}[theorem]{Proposition}

\newtheorem{remark}[theorem]{Remark}

\newcommand \nc{\newcommand}

\nc{\ba}{\begin{array}}\nc{\ea}{\end{array}}
\nc{\be}{\begin{eqnarray}}\nc{\ee}{\end{eqnarray}}
\nc{\beq}{\begin{equation}}\nc{\eeq}{\end{equation}}
\nc{\bex}{\begin{eqnarray*}}\nc{\eex}{\end{eqnarray*}}
\nc{\btm}{\begin{theorem}} \nc{\etm}{\end{theorem}}
\nc{\blm}{\begin{lemma}} \nc{\elm}{\end{lemma}}
\nc{\R}{\mathbb{R}}  \nc{\ld}{\lambda}
\nc{\va}{\varphi}
\nc{\ve}{\varepsilon}

\def\pf{\noindent{\bf Proof.\quad}}\def\endpf{\hfill$\Box$}
\def\di{\mbox{div\,}}

\begin{document}

\title{Stationary solutions with vacuum for a hyperbolic-parabolic chemotaxis model in dimension two}

\author{
Sophia Hertrich\footnote{Department of Mathematics, Wayne State University, Detroit, MI 48202, U.S.A.
%Email: sophia.hertrich@wayne.edu
}
\quad 
Tao Huang$^*$
%\footnote{Department of Mathematics, Wayne State University, Detroit, MI 48202, U.S.A. Email: taohuang@wayne.edu} 
\quad 
Diego Y\'epez
\footnote{Center for Communication Research, Princeton, NJ 08540, U.S.A.}
\footnote{
Department of Mathematics, Princeton University, Princeton, NJ 08540, U.S.A.}
\quad 
Kun Zhao
\footnote{College of Mathematical Sciences, Harbin Engineering University, Harbin, Heilongjiang 150001, P.R. China}
\footnote{Department of Mathematics, Tulane University, New Orleans, LA 70118, U.S.A. 
%Email: 
}
%\date{}
}
\maketitle

\begin{abstract} 

In this research, we study the existence of stationary solutions with vacuum to a hyperbolic-parabolic chemotaxis model with nonlinear pressure in dimension two that describes vasculogenesis. We seek solutions in the radial symmetric class of the whole space, in which the system will be reduced to a system of ODE’s on $(0,\infty)$. The fundamental solutions to the ODE system are the Bessel functions of different types. We find two nontrivial solutions. One is formed by half bump (positive density region) starting at $r=0$ and a region of vacuum on the right. Another one is a full nonsymmetric bump away from $r=0$. These solutions bear certain resemblance to {\it in vitro} vascular network and the numerically produced structure by Gamba et al \cite{Gamba03}. We also show the nonexistence of full bump starting at $r=0$ and nonexistence of full symmetric bump away from $r=0$. 
\end{abstract}

\tableofcontents

\section{Introduction}
\subsection{History}

During the biological process of chemotaxis, cells migrate based on the gradient of a chemical in the environment. In positive chemotaxis, cells move toward the chemical called a chemoattractant. Chemotaxis is basic for life and versatile in function as cells, both prokaryotic and eukaryotic, must be able to react and move based on external stimuli \cite{chemotaxis}. Consider, for example, a white blood cell moving toward an injury site, sperm migrating toward an egg, and bacteria moving toward a food source \cite{chemotaxis}. Furthermore, positive chemotaxis is a key mechanism in various morphogenic events such as our focus, vasculogenesis \cite{Roman00}.
During vasculogenesis, endothelial cells secrete chemoattractants that allow endothelial cells to cluster, thus organizing into vascular networks \cite{CHS23}. 

The following hyperbolic-parabolic chemotaxis model with nonlinear pressure has been proposed in \cite{Gamba03} to describe vasculogenesis
\beq\label{flow2D}
\left\{
\begin{array}{ll}
\rho_t+\di(\rho u)=0\\
\\
(\rho u)_t+\di(\rho u\otimes u)+\nabla p(\rho)=\chi\rho\nabla \phi-\alpha \rho u\\
\\
\delta \phi_t=D\Delta\phi+a\rho-b\phi.
\end{array}
\right.
\eeq
It reproduces key features of experiments of in vitro formation of blood vessels showing that cells randomly spread on a gel matrix autonomously organize to a connected vascular network \cite{Serini03}. 
Here the coefficients satisfy $D,\chi, \alpha>0$, $a,b,\delta\geq 0$, and moreover
$\rho(x,y,t):\R^2\times(0,\infty)\rightarrow \R$ denotes the density of endothelial cells, 
$u(x,y,t): \R^2\times(0,\infty)\rightarrow \R^2$ denotes the velocity of endothelial cells, 
and $\phi(x,y,t):\R^2\times(0,\infty)\rightarrow \R$ denotes the concentration of the chemoattractant secreted by the endothelial cells. The nonlinear pressure
\beq\label{defP}
p(\rho)=\frac{\ve'}{2}\rho^{\gamma}
=\frac{\ve(\gamma-1)}{\gamma}\rho^{\gamma}
\eeq 
for $\gamma>1$ and $\ve,\ve'>0$, is related to the fact that closely packed cells resist compression due to the impenetrability of cellular matter. The term $\alpha\rho u$ corresponds to a damping force and the coefficient $\alpha>0$ is a result of the interaction between cells and the underlying substratum. The quantity $\chi>0$ indicates that the cell is attracted along the chemical concentration gradient. The constants $\delta$ and $D$ are the relaxation time scale and the diffusion coefficient of the chemoattractant, respectively. Furthermore, the nonnegative
constants $a$ and $b$ denote the death and secretion rates of the chemoattractant, respectively.

Most of the analytical research of system \eqref{flow2D} focuses on studying the global wellposedness and long-time asymptotic behavior of the solutions. In \cite{Russo12, Russo13}, global strong solutions without vacuum have been constructed provided the initial data are {small perturbations of small enough constant ground states, and the solutions are shown to converge to the ground states} as $t\rightarrow \infty$. In \cite{ LW24, LPW221}, the smallness condition on the ground density has been removed and the decay rates have been improved when the pressure function satisfies suitable conditions. More well-posedness results in various spaces have been studied in \cite{CHS23, LPW222}. 

\bigskip
$\textbf{Acknowledgments}$. The third author is supported by the IDA Postdoctoral Fellowship. The third author would like to thank Nick Katz and the Mathematics Department at Princeton University for their hospitality.

%\bigskip
\subsection{Results}
In this paper,  we construct stationary solutions to system \eqref{flow2D} with vacuum on $\R^2$, i.e.,
\beq\label{steqn2D}
\left\{
\begin{array}{ll}
\di(\rho u)=0\\
\\
\di(\rho u\otimes u)+\nabla p(\rho)=\chi\rho\nabla \phi-\alpha \rho u\\
\\
D\Delta\phi+a\rho-b\phi=0.
\end{array}
\right.
\eeq
In dimension one, stationary solutions with vacuum (bump solutions) over $\R$ and in a bounded interval with zero-flux boundary condition have been constructed in \cite{2016Paper, CCWWZ20}. In \cite{2016Paper}, {the authors} also proved that there is one ``single bump" stationary solution that is unique up to translation on $\R$. 
The stability of transition layer solutions of \eqref{flow2D} on $[0,\infty)$ has been proved in 
\cite{HPWZ21}. When $\gamma=2$ in the pressure function $p(\rho)$, similar result has been proved for bounded intervals in \cite{PZ23}.

%We may consider the system on $\R^2$ or on a bounded disk centered at origin $B_R\subset\R^2$ for $0<R<\infty$ with boundary conditions on $\partial B_R$
%\beq\label{bdycond1}
%\frac{D\rho}{D\nu}=\frac{D\phi}{D\nu}=0,\quad u=0,
%\eeq
%where $\nu$ is the unit outer normal vector on $\partial B_R$.

In our paper, we aim to investigate the existence of stationary solutions and construct stationary solutions to a two-dimensional positive chemotaxis model.
To begin, we impose the following assumptions on the stationary model \eqref{steqn2D} on $\R^2$ (see Sec \ref{sec2} for more details)
\begin{itemize}
\item[(i)] We let $\gamma=2$. 
 \item[(ii)] The solutions $\rho, \phi$ {depend only} on the radius, and the velocity field $u$ is axisymmetric. 
\end{itemize}
As a result, the velocity field $u$ vanishes, and the stationary model \eqref{steqn2D} reduces to the following system of
$(\rho(r), \phi(r))$
\beq\label{steqn}
\left\{
\begin{array}{ll}
\partial_r( p(\rho))=\chi\rho\phi_r\\
\\
\displaystyle D\phi_{rr}+D\frac{\phi_r}{r}+a\rho-b\phi=0.
\end{array}
\right.
\eeq
In this paper, we investigate the {existence of} nontrivial solutions to this system {satisfying $\phi \geq \rho\geq 0$} with suitable boundary values on $[0,\infty)$ with the following regularity, 
\beq\label{regcond}
\rho\in C^0 [0,\infty),\quad \phi\in C^2 [0,\infty).
%previous was (
\eeq

%USE FOR RESEARCH, NOT THIS PAPER
%with boundary conditions (when the domain is bounded)
%\beq
%\rho(0)=\rho_0, \quad \phi(0)=\phi_0,\quad %\rho_r(R)=\phi_r(R)=0.
%\eeq

To understand the system biologically, we recall $\rho$ denotes the density of endothelial cells and $\phi$ denotes the concentration of the chemoattractant secreted by the endothelial cells. $D>0$ is the diffusion coefficient of the chemoattractant. The quantity $\chi>0$ indicates that the cell is attracted along the chemical concentration gradient. The nonnegative parameters $a$ and $b$ are the death and secretion rates of the chemoattractant respectively. The function $p(\rho)$ is given by the pressure law for isentropic gases, defined as follows for $\gamma=2$,
%info from 148 - FAR ENOUGH OFF?
\beq\label{defP2}
p(\rho)=\frac{\varepsilon'}{2}\rho^2.
\eeq
For convenience, we introduce $\beta\geq0$ such that 
\beq\label{defbeta}
b=D\beta^2.
\eeq

%In the regions where the cell density $\rho=0$, we further divide investigating the solution depending on a constant $\beta$ which we define to be $\sqrt{\frac{b}{D}}$ for convenience. 

The fundamental solutions to the system \eqref{steqn} are related to the Bessel functions of different types. This makes our constructions more complicated than those in dimension one in \cite{2016Paper}, in which the fundamental solutions are in forms of sine, cosine and exponential functions. 

\section{Derivation of equations}\label{sec2}
\setcounter{equation}{0}

To derive the system \eqref{steqn} from \eqref{steqn2D}, we assume that {the} stationary solutions $\rho$ and $\phi$ depend {only} on the radius $r=\sqrt{x^2+y^2}$. We also require that the velocity field is axisymmetric 
$$
u(r, \theta)=\big(v(r)\cos \theta,\, v(r)\sin\theta\big),
$$
where $v(r): [0,\infty)\rightarrow \R$ is a scalar function. For compatibility, we need to assume 
\beq
v(0)=0.
\eeq 
By the definition of polar coordinates
$$
x=r\cos\theta,\quad y=r\sin\theta
$$
we have
$$
r_x=\frac{x}{r}=\cos\theta,\quad 
r_y=\frac{y}{r}=\sin\theta,\quad
\theta_x=\frac{-y}{x^2+y^2}=-\frac{\sin\theta}{r},\quad
\theta_y=\frac{x}{x^2+y^2}=\frac{\cos\theta}{r}.
$$
Hence
$$
\partial_x=\cos\theta\,\partial_r-\frac{\sin\theta}{r}\,\partial_{\theta},\quad
\partial_y=\sin\theta\,\partial_r+\frac{\cos\theta}{r}\,\partial_{\theta}.
$$
Direct computation implies 
\beq\notag
\di(\rho u)=(\rho v\cos\theta)_x+(\rho v\sin\theta)_y
=(\rho v)_r+\frac{\rho v}{r}=\frac{(r\rho v)_r}{r}=0.
\eeq
Thus, 
\beq\notag
\rho v=\frac{C}{r}
\eeq
for some constant $C$. Since we need to find regular stationary solutions (e.g. assumption \eqref{regcond}) and the boundary value $v(0)=0$, we need to pick $C\equiv 0$. Therefore, $\rho v\equiv 0$, which implies $\rho u\equiv 0$.
Thus, the system \eqref{steqn2D} is reduced to the system \eqref{steqn} of $(\rho(r), \phi(r))$. 

%with boundary conditions (when the domain is bounded)
%\beq\label{bdycond2}
%\rho(0)=\rho_0, \quad \phi(0)=\phi_0,\quad \rho_r(R)=\phi_r(R)=0,
%\eeq
%It is more natural to consider this boundary condition than non flux boundary condition at $r=0$.
%or maybe nonflux boundary conditions on both sides
%\beq\label{bdycond22}
%\rho_r(0)= \phi_r(0)=0,\quad \rho_r(R)=\phi_r(R)=0,
%\eeq

\medskip
Before the end of this section, we state some facts on the energy of smooth solutions to the  system \eqref{steqn}.
By the computation in Appendix \ref{App1}, the energy of the hyperbolic-parabolic chemotaxis system \eqref{flow2D} is given as follows
\beq\notag
E(t)=\int \left( \frac{1}{2}\rho |u|^2+\frac{\ve}{2}\rho^{2}+\frac{\chi D}{2a}|\nabla \phi|^2+\frac{\chi b}{2a}|\phi|^2-\chi \rho\phi\right)\,d{\bf x}.
\eeq
Here, we take $\gamma=2$, and denote ${\bf x}=(x,y)\in \R^2$ and $\int\,d{\bf x}=\int_{\R^2}\,d{\bf x}$ for simplicity. If we only consider stationary solutions to system \eqref{steqn2D}, we can derive the following estimate by multiplying $\frac{\chi}{a}\phi$ to both sides of equation of $\phi$ and integrating by parts over $\R^2$
\beq\notag
\int \left(\frac{\chi D}{a}|\nabla \phi|^2+\frac{\chi b}{a}|\phi|^2\right)\,d{\bf x}
=\int \chi\rho\phi\,d{\bf x}.
\eeq
Therefore, the energy of stationary system \eqref{steqn2D} is 
\beq\label{engst2D}
E_s:=\int \left( \frac{1}{2}\rho |u|^2+\frac{\ve}{2}\rho^{2}+\frac{\chi}{2} \rho\phi\right)\,d{\bf x}.
\eeq
We will compute the energy of the nontrivial solutions to system \eqref{steqn} that we construct below. It will be negative, smaller than the constant stationary solutions when $u\equiv 0$, $\rho$ and $\phi$ are constants.

\section{General solutions} 
\setcounter{equation}{0}

In this section, we first investigate the general solutions to the system \eqref{steqn} on $\R$. By the first equation of \eqref{steqn}, we have
\beq\label{eqnrho0}
p(\rho)_r=\ve\rho\rho_r=\chi\rho\phi_r.
\eeq
We divide it into two different regions. 

\subsection*{When $\rho=0$} 
The equation of $\rho$ becomes trivial and the second equation becomes
\beq\notag
D\phi_{rr}+D\frac{\phi_r}{r}-b\phi=0.
\eeq
By the definition of $\beta\geq 0$ in \eqref{defbeta}, {this equation is} reduced as follows
\beq\label{eqnphi0}
\phi_{rr}+\frac{\phi_r}{r}-\beta^2\phi=0.
\eeq
When $\beta\equiv 0$, the general solution is 
\beq\label{solph00}
\phi(r)=c_1\ln r+c_2
\eeq
 for any constants $c_1$ and $c_2$.
When $\beta>0$, the general solution is 
\beq\label{solph01}
\phi(r)=A_1 J_0(i\beta r)+A_2Y_0(-i\beta r)
\quad
\mbox{or}
\quad
\phi(r)=A_1 I_0(\beta r)+A_2K_0(\beta r).
\eeq
Here $J_0$ is the Bessel function of first kind 
\beq\label{defB1}
J_0(r)=\sum\limits_{k=0}^\infty \frac{(-1)^k}{k!(k+1)!}\left(\frac{r}{2}\right)^{2k} .
\eeq
$Y_0$ is the Bessel function of second kind
\beq\label{defB2}
Y_0(r)=\frac{2}{\pi}\left(\left(\ln\left(\frac{r}{2}+\gamma_{em}\right)\right)J_0(r)+\sum\limits_{k=1}^\infty \frac{(-1)^{k+1}H_k}{(k!)^2}\left(\frac{r}{2}\right)^{2k}\right),
\eeq
where $\gamma_{em}$ is the Euler-Mascheroni constant and $H_k=\sum\limits_{l=1}^{k}\frac{1}{l}$ is a harmonic number. 
$I_0$ is the modified Bessel function of first kind
\beq\label{defMB1}
I_0(r)=\sum\limits_{k=0}^\infty \frac{1}{(k!)^2}\left(\frac{r}{2}\right)^{2k}.
\eeq
$K_0$ is the modified Bessel function of second kind 
\beq\label{defMB2}
K_0(r)=\int_0^\infty\cos(r\sinh t)\,dt
=\int_0^\infty\frac{\cos(rt)}{\sqrt{t^2+1}}\,dt.
\eeq

\subsection*{When $\rho\neq 0$}
The equation \eqref{eqnrho0} can be written as follows
\beq\notag
\ve\rho_r=\chi\phi_r.
\eeq
Integrating both sides, we have 
\beq\label{eqnK}
\ve\rho=\chi\phi+K
\eeq
for some constant $K$ to be determined later. This implies
\beq\label{RPK}
\rho=\frac{\chi}{\ve}\phi+\frac{K}{\ve}.
\eeq
Plugging into the second equation of \eqref{steqn}, we have 
\beq\notag
D\phi_{rr}+D\frac{\phi_r}{r}+a\frac{\chi\phi+K}{\ve}-b\phi
=D\phi_{rr}+D\frac{\phi_r}{r}+\left(\frac{a\chi}{\ve}-b\right)\phi+\frac{aK}{\ve}
=0.
\eeq
By the assumption $b=D\beta^2$, it holds
\beq\label{eqnphi1}
\phi_{rr}+\frac{\phi_r}{r}+\left(\frac{a\chi}{D\ve}-\beta^2\right)\phi+\frac{aK}{D\ve}
=0.
\eeq
The fundamental solutions to the ODE \eqref{eqnphi1} are different in the following three cases. 
\medskip

\noindent{\bf Case 1.} When $\frac{a\chi}{D\ve}-\beta^2=0$, the general solution of \eqref{eqnphi1} is 
\beq\label{solph11}
\phi(r)=c_1\ln r+c_2-\frac{aK}{4D\ve}r^2.
\eeq
Thus, by the equation \eqref{RPK}, we obtain
\beq\label{solrho11}
\rho(r)=\frac{c_1\chi}{\ve}\ln r+\frac{c_2\chi}{\ve}-\frac{\chi aK}{4D\ve^2}r^2+\frac{K}{\ve}.
\eeq

\medskip
\noindent{\bf Case 2.} When $\frac{a\chi}{D\ve}-\beta^2<0$, for convenience, we assume that 
$$
\frac{a\chi}{D\ve}-\beta^2=-\xi^2
$$
for some $\xi>0$.
Similar to \eqref{solph01}, the general solution is 
\beq\label{solph12}
\phi(r)=c_1 I_0(\xi r)+c_2K_0(\xi r)+\frac{aK}{D\ve\xi^2},
\eeq
\beq\label{solroh12}
\rho(r)=\frac{c_1\chi}{\ve}I_0(\xi r)+\frac{c_2\chi}{\ve}K_0(\xi r)+\frac{\chi aK}{D\ve^2\xi^2}+\frac{K}{\ve}.
\eeq

\medskip
\noindent{\bf Case 3.} When $\frac{a\chi}{D\ve}-\beta^2>0$, for convenience, we assume that 
$$
\frac{a\chi}{D\ve}-\beta^2=\omega^2
$$
for some $\omega>0$.
The general solution is 
\beq\label{solph13}
\phi(r)=c_1 J_0(\omega r)+c_2Y_0(\omega r)-\frac{aK}{D\ve\omega^2},
\eeq
\beq\label{solpho13}
\rho(r)=\frac{c_1\chi}{\ve}J_0(\omega r)+\frac{c_2\chi}{\ve}Y_0(\omega r)-\frac{\chi aK}{D\ve^2\omega^2}+\frac{K}{\ve}.
\eeq

\section{Transitions between vacuum and positivity regions} 
\setcounter{equation}{0}

This section {is} devoted to investigate the transitions between vacuum and positivity regions of density. Our main result is as follows.

\begin{proposition}\label{proptran} For any solutions $\rho\in C^0$ and $\phi\in C^2$ of system \eqref{steqn}, suppose that there is transition point $\bar r>0$ between vacuum and positivity regions of density $\rho(r)$.  Then, the solution $\phi(r)$ is $C^2$ at $\bar r$ if and only if 
 $\phi(r)$ is $C^1$ at $\bar r$ and 
 \beq\label{nescod11}
\beta>0,\quad \phi(\bar r)=-\frac{K}{\chi}\geq 0.%=\phi_0-\frac{\ve}{\chi}\rho_0.
\eeq
\end{proposition}

\begin{remark}
(1) From now on, we only assume that $\beta>0$.

\noindent (2) Since $\phi (r)\geq 0$ for $r\geq 0$, the equation \eqref{nescod1} below implies that we require
\beq\label{Kneg}
K\leq 0.
\eeq

\noindent (3) The case $\bar r=0$ will be considered later.
\end{remark}

\pf We first consider the case $\beta=0$. On the vacuum side, by \eqref{solph00}, we have $\phi(r)\equiv c$ for some constant $c\geq 0$. On the non-vacuum side, we have two different subcases. First when $a=0$, we obtain by \eqref{solph11} 
\beq\notag
\phi(r)=c_1 \ln r+c_2.
\eeq
The only possible case for which $\phi \in C^2$ is $\phi\equiv c$ on the non-vacuum side. Then $\rho(r)=\frac{c\chi}{\ve}+\frac{K}{\ve}$ which never vanishes. When $a>0$, if we want to make $\phi \in C^2$ at $\bar r$, we {must} have $\phi(r)\equiv c$ by the equation \eqref{eqnphi1}. This implies $\rho(r)=\frac{c\chi}{\ve}+\frac{K}{\ve}$ which never vanishes. Therefore, we need to assume that $\beta>0$ if we want to make $\phi\in C^2$.

When $\beta>0$ and if we assume that $\phi\in C^2$ at $\bar r$, we may compute the derivatives at $\bar r>0$ on the vacuum side as follows (for simplicity, we omit the notation of right-side derivative)
\beq\label{eqnA1}
\phi(\bar r)
=A_1I_0(\beta  \bar r)+A_2K_0(\beta  \bar r),
\eeq
\beq\label{eqnA2}
\phi'( \bar r)=A_1\partial_r I_0(\beta  \bar r)+A_2 \partial_r K_0(\beta  \bar r),
\eeq
\beq\label{eqnA3}
 \phi''( \bar r)=A_1\partial_{rr} I_0(\beta  \bar r)+A_2\partial_{rr} K_0(\beta  \bar r).
\eeq
The linear system of $A_1$ and $A_2$ has at least one solution only when the following equation is valid on the vacuum side
\beq\label{cdphi}
 \bar r^2\phi''( \bar r)+ \bar r\phi'( \bar r)-\beta^2 \bar r^2\phi( \bar r)=0.
\eeq
Here we have used that $I_0$ and $K_0$ are fundamental solutions to the equation \eqref{eqnphi0}. By the assumption $\phi\in C^2$ and the equation \eqref{eqnphi0}, we see that \eqref{cdphi} is always valid on the vacuum side.

Let {us} now investigate the condition \eqref{cdphi} on the non-vacuum side. In fact, by the equation \eqref{RPK} and $\rho(\bar r)=0$, we have 
\beq\label{nescod1}
\phi( \bar r)=-\frac{K}{\chi},%=\phi_0-\frac{\ve}{\chi}\rho_0.
\eeq
which implies 
\beq\notag
\frac{a\chi}{D\ve}\phi(\bar r)+\frac{aK}{D\ve}=0.
\eeq
On the non-vacuum side, it holds by equation \eqref{eqnphi1}
\beq\notag
r^2\phi_{rr}+r\phi_r+r^2\left(\frac{a\chi}{D\ve}-\beta^2\right)\phi+\frac{aK}{D\ve}r^2=0.
\eeq
Combining these two equations, we conclude that the condition \eqref{cdphi} is always valid at the transition point $\bar r>0$ (on both sides).

Inversely, if we assume $\phi$ is $C^1$ at $\bar r$ and \eqref{nescod1}, then by the equations \eqref{eqnphi0} and \eqref{eqnphi1} on both sides of $\bar r$, it is not hard to see that $\phi$ is $C^2$ at $\bar r$. 

\endpf

\bigskip
We can solve for $A_1$ and $A_2$ in \eqref{eqnA1}-\eqref{eqnA3} by the first two equations and the Cramer rule
\beq\label{solAs}
A_1=\frac{W_1( \bar  r)}{W(  \bar  r)},\quad A_2=\frac{W_2(  \bar  r)}{W( \bar  r)}
\eeq
where 
\beq\label{defW}
W( r)=
\left|
\begin{array}{cc}
I_0(\beta r)&K_0(\beta r)\\
\partial_rI_0(\beta r)&\partial_r K_0(\beta r)
\end{array}
\right|
\eeq
is the Wronskian of $I_0$ and $K_0$ at $\beta r$, which is not zero since they are fundamental solutions to a linear ODE, and 
\beq\label{W1W2}
W_1(r)=
\left|
\begin{array}{cc}
\phi( r)&K_0(\beta r)\\
\phi'(r)&\partial_r K_0(\beta r)
\end{array}
\right|
\quad
W_2(r)=
\left|
\begin{array}{cc}
I_0(\beta r)&\phi( r)\\
\partial_rI_0(\beta r)&\phi'( r)
\end{array}
\right|.
\eeq
Since $I_0(r)>0$ is increasing and $K_0(r)> 0$ is decreasing for any $0< r<\infty$,  we have 
\beq\notag
W(r)<0.
\eeq

\section{Half bump solutions at \texorpdfstring{$r = 0$}{Lg}} 
\setcounter{equation}{0}

In the rest of {this} paper, we construct solutions to system \eqref{steqn} with vacuum. We first consider a half bump at $r=0$. More precisely, we assume that  $\rho(0)=\rho_0> 0$, $\phi(0)=\phi_0> 0$. %and {\color{red}$\phi_0\geq \rho_0$ ???}. 
There exists some $r_0>0$ such that $\rho> 0$ on $[0,r_0)$ and $\rho\equiv 0$ on $[r_0,\infty)$.
By the assumption of the initial data and \eqref{RPK}, it is not hard to derive
\beq\label{defK0}
K=\ve\rho_0-\chi\phi_0.
\eeq
We also need to check the compatibility of all constants.
First the condition \eqref{Kneg} requires
\beq\label{KNinitial}
 \phi_0\geq \frac{\ve}{\chi}\rho_0.
\eeq
Since when $\rho>0$, we have 
\beq\notag
\rho=\frac{\chi}{\ve}\phi+\frac{K}{\ve}
\eeq
which implies 
\beq\label{Pnescon1}
\phi\geq -\frac{K}{\chi}=\phi_0-\frac{\ve}{\chi}\rho_0\geq 0,
\eeq
which is compatible.

%{\color{red}The background of the model implies that we  require $\phi\geq \rho$ ???}. This in turn implies 
%\beq\notag
%\left(1-\frac{\chi}{\ve}\right)\phi\geq \frac{K}{\ve}.
%\eeq
%When $\frac{\rho_0}{\phi_0}\leq \frac{\chi}{\ve}\leq 1$, this is trivial. When $\frac{\chi}{\ve}> 1$, we require
%\beq\notag
%\phi\leq \frac{-K}{\chi-\ve}.
%\eeq
%It is easy to see that $\frac{-K}{\chi-\ve}\geq-\frac{K}{\chi}$ since $\chi>\ve$. This is also compatible. 

Now, we are ready to construct the nontrivial half bump solutions in three different cases.
\smallskip

\noindent{\bf Case 1.} When $\frac{a\chi}{D\ve}-\beta^2=0$,   the solution $\rho$ is given as follows on $[0,r_0]$ 
\beq\notag
\rho(r)=\frac{c_1\chi}{\ve}\ln r+\frac{c_2\chi}{\ve}-\frac{\chi aK}{4D\ve^2}r^2+\frac{K}{\ve}.
\eeq
Since $\rho_0$ is finite, it is easy to see that 
\beq\notag
c_1=0,\quad c_2=\frac{\ve}{\chi}\left(\rho_0-\frac{K}{\ve}\right)=\phi_0.
\eeq 
Thus, on $[0,r_0]$
\beq\notag
\rho(r)=\rho_0-\frac{\chi aK}{4D\ve^2}r^2,\quad 
\phi(r)=\phi_0-\frac{aK}{4D\ve}r^2.
\eeq
By the necessary condition \eqref{KNinitial}, we have $K\leq 0$, which implies $\rho$ never equals zero. Therefore, there is no half bump solution in this case.

\bigskip
\noindent{\bf Case 2.} When $\frac{a\chi}{D\ve}-\beta^2=-\xi^2<0$, 
the solutions $\rho$ and $\phi$ are given as follows on $[0, r_0]$
\beq\notag
\rho(r)=\frac{c_1\chi}{\ve}I_0(\xi r)+\frac{c_2\chi}{\ve}K_0(\xi r)+\frac{\chi aK}{D\ve^2\xi^2}+\frac{K}{\ve},
\eeq
\beq\notag
\phi(r)=c_1 I_0(\xi r)+c_2K_0(\xi r)+\frac{aK}{D\ve\xi^2}.
\eeq
Since $\rho_0$ is finite, we must let 
%K_0 diverges at 0 so c_2=0 and c_1 comes from setting \rho(0)=\rho_0
\beq\notag
c_2=0.
\eeq 
By observing $\rho_0$ and $\phi_0$, we find
\beq\notag
c_1=\frac{\ve}{\chi}\left(\rho_0-\frac{\chi aK}{D\ve^2\xi^2}-\frac{K}{\ve}\right)=\phi_0-\frac{aK}{D\ve\xi^2}>0
\eeq
since $K\leq 0$ by \eqref{KNinitial}.
Thus, on $[0, r_0]$
\beq\notag
\rho(r)=\left(\rho_0-\frac{\chi aK}{D\ve^2\xi^2}-\frac{K}{\ve}\right)I_0(\xi r)+\frac{\chi aK}{D\ve^2\xi^2}+\frac{K}{\ve}, 
\eeq
\beq\notag
\phi(r)=\left(\phi_0-\frac{aK}{D\ve\xi^2}\right)I_0(\xi r)+\frac{aK}{D\ve\xi^2}.
\eeq
Since $\rho_0>0$ and $I(\xi r)$ is increasing with respect to $r$, we conclude that $\rho$ will never be zero. Therefore, there is no half bump solution in this case.

\bigskip
\noindent{\bf Case 3.} When $\frac{a\chi}{D\ve}-\beta^2=\omega^2>0$, the solutions $\rho$ and $\phi$ are given as follows on $[0, r_0]$
\beq\notag
\rho(r)=\frac{c_1\chi}{\ve}J_0(\omega r)+\frac{c_2\chi}{\ve}Y_0(\omega r)-\frac{\chi aK}{D\ve^2\omega^2}+\frac{K}{\ve},
\eeq
\beq\notag
\phi(r)=c_1 J_0(\omega r)+c_2Y_0(\omega r)-\frac{aK}{D\ve\omega^2}.
\eeq
Since $\rho_0$ is finite, we must let 
\beq\notag
c_2=0.
\eeq 
By observing $\rho_0$ and $\phi_0$, we find
\beq\notag
c_1=\frac{\ve}{\chi}\left(\rho_0+\frac{\chi aK}{D\ve^2\omega^2}-\frac{K}{\ve}\right)=\phi_0+\frac{aK}{D\ve\omega^2}.
\eeq
Thus, on $[0, r_0]$
\beq\notag
\rho(r)=\left(\rho_0+\frac{\chi aK}{D\ve^2\omega^2}-\frac{K}{\ve}\right)J_0(\omega r)-\frac{\chi aK}{D\ve^2\omega^2}+\frac{K}{\ve}, 
\eeq
\beq\notag
\phi(r)=\left(\phi_0+\frac{aK}{D\ve\omega^2}\right)J_0(\omega r)-\frac{aK}{D\ve\omega^2}.
\eeq
Since $J_0$ is decreasing around $r=0$, we require $c_1>0$ in order to obtain a zero point of $\rho$. This is equivalent to 
\beq\label{excod1}
\phi_0\geq -\frac{aK}{D\ve\omega^2}=-\left(\frac{\beta^2}{\chi\omega^2}+\frac{1}{\chi}\right)K.
\eeq
This is compatible with \eqref{KNinitial} and \eqref{Pnescon1}. 
Direct computation implies 
\beq\notag
L:=-\frac{\chi aK}{D\ve^2\omega^2}+\frac{K}{\ve}
=-\frac{K}{\ve}\left(\frac{\chi a}{D\ve\omega^2}-1\right)=-\frac{K}{\ve}\frac{\beta^2}{\omega^2}\geq0.
\eeq
By the fact that $c_1>0$, it holds
\beq\notag
\rho_0>L.
\eeq
Suppose that $J_0(\omega r)$ {obtains} its first minimum value $-m$ at $\tilde r_0>0$ for some positive $m$. To guarantee there is a zero of $\rho$ in some interval $0<r_0\leq \tilde r_0$, we need 
\beq\notag
\frac{L}{\rho_0-L}\leq m
\eeq
which is equivalent to 
\beq\notag
-\frac{K}{\ve}\frac{\beta^2}{\omega^2}\leq \frac{m\rho_0}{1+m}
\eeq
or 
\beq\label{excon2}
\frac{\chi}{\ve}\frac{\beta^2}{\omega^2}\phi_0\leq \left(\frac{m}{1+m}+\frac{\beta^2}{\omega^2}\right)\rho_0.
\eeq
This is compatible with \eqref{KNinitial}. The zero point $r_0$ can be computed as follows
\beq\notag
J(\omega r_0)=-\frac{L}{\rho_0-L}.
\eeq

Now, we need to make $\phi\in C^2$ at $r_0$. By the discussion in last several sections, when $\rho=0$, it holds
\beq\notag
\phi=A_1 I_0(\beta r)+A_2K_0(\beta r).
\eeq
The solution $\phi\in C^2$ at $r_0$ when $A_1$ and $A_2$ are given by 
\beq\label{solAsr0}
A_1=\frac{W_1(r_0)}{W(r_0)},\quad A_2=\frac{W_2(r_0)}{W(r_0)}
\eeq
where $W(r_0)$, $W_1(r_0)$ and $W_2(r_0)$ are given by \eqref{defW} and \eqref{W1W2}. Since $I_0(r)>0$ is increasing and $K_0(r)> 0$ is decreasing for any $0< r<\infty$,  we have 
\beq\notag
W(r_0)<0.
\eeq
Since $\phi>0$ is decreasing around $r_0$, we have 
\beq\notag
W_2(r_0)<0,
\eeq
which implies 
\beq\label{A2pos}
A_2>0.
\eeq
Since $I_0$ approaches $\infty$ when $r\rightarrow\infty$, we need to make $A_1=0$, which is equivalent to 
\beq\label{5.34}
\begin{split}
&0=W_1( r_0)=
\left|
\begin{array}{cc}
\phi(r_0)&K_0(\beta r_0)\\
\phi'( r_0)&\partial_r K_0(\beta r_0)
\end{array}
\right|
=\phi(r_0)\partial_r K_0(\beta r_0)-\phi'(r_0)K_0(\beta r_0).
%=&\left(\left(\phi_0+\frac{aK}{D\ve\omega^2}\right)J_0(\omega r_0)-\frac{aK}{D\ve\omega^2}\right)\partial_r K_0(\beta r_0)-\left(\phi_0+\frac{aK}{D\ve\omega^2}\right)\partial_r J_0(\omega r_0)K_0(\beta r_0)\\
%=&-\frac{K}{\chi}\partial_r K_0(\beta r_0)-\left(\phi_0+\frac{aK}{D\ve\omega^2}\right)\partial_r J_0(\omega r_0)K_0(\beta r_0).
\end{split}
\eeq
%Here we have used the fact \eqref{nescod1} in last identity.
Since $A_1=0$, to make $\phi$ being $C^2$ at $r_0$, we obtain by Proposition \ref{proptran},
\beq\label{5.35}
\phi(r_0)=A_2K_0(\beta r_0)=-\frac{K}{\chi},
\eeq
\beq\label{5.36}
\phi'(r_0)=A_2\partial_rK_0(\beta r_0)=\left(\phi_0+\frac{aK}{D\ve\omega}\right)\partial_r J_0(\omega r_0).
\eeq
Observe that the equation \eqref{5.36} is valid if we have  \eqref{5.34} and \eqref{5.35}. Direct computation implies \eqref{5.35} can be written as follows
\beq\label{5.37}
\begin{split}
&W_2(r_0)K_0(\beta r_0)
=K_0(\beta r_0)\left(c_1I_0(\beta r_0)\partial_r J_0(\omega r_0)+\frac{K}{\chi}\partial_rI_0(\beta r_0)\right)\\
=&-\frac{K}{\chi} W(r_0)
=-\frac{K}{\chi}\left(I_0(\beta r_0)\partial_r K_0(\beta r_0)-\partial_rI_0(\beta r_0)K_0(\beta r_0)\right).
\end{split}
\eeq
It is not hard to see that \eqref{5.34} is equivalent to
\beq\notag
-\frac{K}{\chi}\partial_r K_0(\beta r_0)=
%\left(\phi_0+\frac{aK}{D\ve\omega}\right)
c_1\partial_r J_0(\omega r_0)K_0(\beta r_0),
\eeq
which implies \eqref{5.37} is always valid.
%Tao had next sentence in red
Therefore, we should have the following requirement at $r_0$ to make $\phi \in C^2$
\beq\notag
-\frac{K}{\chi}\partial_r K_0(\beta r_0)=
%\left(\phi_0+\frac{aK}{D\ve\omega}\right)
c_1\partial_r J_0(\omega r_0)K_0(\beta r_0),
\eeq
which is still possible by the definition of Bessel functions. 
Therefore, we obtain a nontrivial half bump solution to system \eqref{steqn} in this case. 

Now, let us compute the energy of this solution. Since $u\equiv 0$, the energy $E_s$ in \eqref{engst2D} of the system becomes
\beq\notag
\int \left( \frac{\ve}{2}\rho^{2}+\frac{\chi}{2} \rho\phi\right)\,d{\bf x}
=\int \frac{\rho}{2}\left( \ve\rho+\chi \phi\right)\,d{\bf x}.
\eeq
For simplicity, we use Cartesian coordinates instead of polar coordinates. By the definition of $K$ in \eqref{eqnK}, we obtain
\beq\notag
\int \left( \frac{\ve}{2}\rho^{2}+\frac{\chi}{2} \rho\phi\right)\,d{\bf x}
=\int \frac{1}{2}\rho K\,d{\bf x}.
\eeq
By the construction above, it is not hard to check that $K$ is negative. Otherwise, it is impossible to make $\phi\in C^2$ at the transition point between vacuum and positivity regions. Therefore, the energy of the constructed solution is always negative.  

\section{Single whole bump solutions touching \texorpdfstring{$r = 0$}{Lg}} 
\setcounter{equation}{0}

In this section, we show that there is no single whole bump solutions touching $r=0$ to the system \eqref{steqn}. More precisely, we assume $\rho(0)=\rho_0=0$ and $\phi(0)=\phi_0\geq 0$. There exists some $r_0>0$ such that $\rho> 0$ on $(0,r_0)$ and $\rho\equiv 0$ on $[r_0,\infty)$. We need to investigate in three different cases. 

\bigskip
\noindent{\bf Case 1.} When $\frac{a\chi}{D\ve}-\beta^2=0$, 
the solutions $\rho$ and $\phi$ are given as follows on $[0, r_0]$, 
\beq\notag
\rho(r)=\frac{c_1\chi}{\ve}\ln r+\frac{c_2\chi}{\ve}-\frac{\chi aK}{4D\ve^2}r^2+\frac{K}{\ve},
\eeq
\beq\notag
\phi(r)=c_1\ln r+c_2-\frac{aK}{4D\ve}r^2.
\eeq 
Since $\rho_0=0$, we must have 
\beq\notag
c_1=0, %to elimate ln(r) term
\quad
c_2=-\frac{K}{\chi}. %set \rho =0 and solve for c_2
\eeq
Thus, on $[0, r_0]$, we have 
\beq\notag
\rho(r)=-\frac{\chi aK}{4D\ve^2}r^2.
\eeq
%Should be good, just double check depneding on reg condition
Notice, $\rho(r)=0$ only when $r=0$. Thus, we have no whole bump solution in this case.
%Conditions for r_0: r_0>0 st \rho(r_0)=0

\bigskip
\noindent{\bf Case 2.} When $\frac{a\chi}{D\ve}-\beta^2=-\xi^2<0$, 
the solutions $\rho$ and $\phi$ are given as follows on $[0, r_0]$

\beq\notag
\rho(r)=\frac{c_1\chi}{\ve}I_0(\xi r)+\frac{c_2\chi}{\ve}K_0(\xi r)+\frac{\chi aK}{D\ve^2\xi^2}+\frac{K}{\ve},
\eeq
\beq\notag
\phi(r)=c_1 I_0(\xi r)+c_2K_0(\xi r)+\frac{aK}{D\ve\xi^2}.
\eeq 
Since $\rho_0=0$, we have
\beq\notag
c_2=0, \quad
%to eliminate K_0 term since it diverges at x=0
c_1=-\frac{\ve}{\chi} \left( \frac{\chi aK}{D\ve^2\xi^2}+\frac{K}{\ve}\right).
%set p_0=0 and solve for c_1
\eeq
Thus, on $[0, r_0]$, 
\beq\notag
\rho(r)= \left(\frac{\chi aK}{D\ve^2\xi^2}+\frac{K}{\ve}\right)\left(1-I_0(\xi r) \right).
\eeq
Notice that $\rho(r)=0$ for $r$ such that $I_0(\xi r)=1$. This only holds when $r=0$, thus, we have no whole bump solution in this case.

\bigskip
\noindent{\bf Case 3.} When $\frac{a\chi}{D\ve}-\beta^2=\omega^2>0$,
the solutions $\rho$ and $\phi$ are given as follows on $[0, r_0]$, 
\beq\notag
\rho(r)=\frac{c_1\chi}{\ve}J_0(\omega r)+\frac{c_2\chi}{\ve}Y_0(\omega r)-\frac{\chi aK}{D\ve^2\omega^2}+\frac{K}{\ve},
\eeq
\beq\notag
\phi(r)=c_1 J_0(\omega r)+c_2Y_0(\omega r)-\frac{aK}{D\ve\omega^2}.
\eeq
Since $\rho_0=0$,
\beq\notag 
c_2=0,
\quad 
c_1=-\frac{\ve}{\chi}\left(-\frac{\chi aK}{D\ve^2\omega^2}+\frac{K}{\ve}\right).
\eeq
Thus, on $[0, r_0],$
\beq\notag
\rho(r)=\left(-\frac{\chi aK}{D\ve^2\omega^2}+\frac{K}{\ve}\right)\left(1-J_0(\omega r)\right).
\eeq
Notice that $\rho(r)=0$ for $r$ such that $J_0(\omega r)=1$. This only holds when $r=0$, thus, we have no whole bump solution in this case.

\section{Single whole bump solutions in \texorpdfstring{$(0,\infty)$}{Lg}} 
\setcounter{equation}{0}

In this section, we construct nontrivial solutions to system \eqref{steqn} with single whole bump inside $(0,\infty)$. More precisely, we assume $\rho(0)=\rho_0=0$ and $\phi(0)=\phi_0\geq 0$. There exists some $r_1>r_0>0$ such that $\rho> 0$ on $(r_0,r_1)$ and $\rho\equiv 0$ otherwise. 

We first show that there is no symmetric solution $\rho$ on $(r_0,r_1)$ about $r=(r_0+r_1)/2$. By the relation \eqref{RPK}, it is not hard to see that $\phi$ is also symmetric on $(r_0,r_1)$, which implies 
\beq\label{eq7.1}
\phi(r_0)=\phi(r_1)=-\frac{K}{\chi},\quad \phi_r(r_0)=-\phi_r(r_1),\quad\phi_{rr}(r_0)=\phi_{rr}(r_1).
\eeq
By the fact $\phi\in C^2$ and equation \eqref{eqnphi0}, we conclude that 
\beq\label{eq7.2}
\phi_r(r_0)=\phi_r(r_1)=0.
\eeq

We first consider the case $\beta\equiv 0$. In this case we have the following general solution on the region of vacuum
$$
\phi(r)=c_1\ln r+c_2.
$$
By the condition \eqref{eq7.2}, we have $c_1\equiv 0$, which implies $\phi(r)=-\frac{K}{\chi}$ on $[0,r_1]\cup[r_1,\infty)$. Since we only consider the solutions with finite energy, we conclude that $\phi(r)\equiv 0$ on $[0,r_1]\cup[r_1,\infty)$, which also implies $K\equiv 0$. Therefore, the equation \eqref{eqnphi1} on the region $\rho\neq 0$ becomes
\beq\notag
\phi_{rr}+\frac{\phi_r}{r}+\frac{a\chi}{D\ve}\phi
=0
\eeq 
with zero boundary conditions on $(r_0, r_1)$. Since the solutions also satisfy the condition \eqref{eq7.2}, it is not hard to see that $\phi(r)\equiv 0$ for any $r\geq 0$, which implies that there is no nontrivial symmetric solution when $\beta\equiv 0$.

\medskip
Now let us consider the case $\beta\neq 0$. In this case we have the following general solution on the region of vacuum
$$
\phi(r)=A_1I_0(\beta r)+A_2K_0(\beta r).
$$
On $[0,r_0]$, we have $A_2\equiv 0$ since $K_0(r)$ is unbounded when $r\rightarrow 0^+$. Since $I_0(r)$ is strictly increasing on $r>0$ and achieves the minimum at $r_0$, we have $A_1\equiv 0$ by the condition \eqref{eq7.2} at $r_0$. This implies $\phi\equiv 0$ on $[0,r_0]$. Similarly, we can also prove that $\phi\equiv 0$ on $[r_1,\infty)$. By the first condition in \eqref{eq7.1}, we conclude that $K=0$ and $\phi(r_1)=0$.
Therefore, the equation \eqref{eqnphi1} on the region $\rho\neq 0$ becomes
\beq\notag
\phi_{rr}+\frac{\phi_r}{r}+\left(\frac{a\chi}{D\ve}-\beta^2\right)\phi
=0
\eeq 
with zero boundary conditions on $(r_0, r_1)$. Since the solutions also satisfy the condition \eqref{eq7.2}, it is not hard to see that $\phi(r)\equiv 0$ for any $r\geq 0$, which implies that there is no nontrivial symmetric solution when $\beta\neq 0$.

The byproduct is the following proposition

\begin{proposition}\label{prop7.1} Suppose that the solutions $\rho, \phi$ satisfy all the assumptions above in this section. Moreover, if we have $K\equiv 0$ and the solution $\phi$ satisfies \eqref{eq7.2}, then there is no 
nontrivial solution to the system \eqref{steqn}.
\end{proposition}

\bigskip
Now, let us seek for non-symmetric solutions on $(r_0,r_1)$. By the definition of $r_0$ and $r_1$, we still have 
\beq\label{eq7.5}
\phi(r_0)=\phi(r_1)=-\frac{K}{\chi}.
\eeq

We first consider the case $\beta\equiv 0$, which implies 
$\phi(r)=c_1\ln r+c_2$ on the vacuum region. Since we seek for finite energy solutions, we must have $c_1\equiv 0$ and $\phi(r)=c_2=-\frac{K}{\chi}\equiv 0$ on $[0,r_1]\cup[r_1,\infty)$. By similar arguments in the symmetric case, we conclude that there is no nontrivial solution when $\beta\equiv 0$.

Now, let us consider the case $\beta\neq 0$. By the discussion in the symmetric case, we know 
\beq\label{eq7.6}
\phi(r)=A_1 I_0(\beta r),\ \mbox{on } [0,r_0],\quad
\phi(r)=A_2 K_0(\beta r), \ \mbox{on } [r_1,\infty).
\eeq
By the initial data, it holds $A_1=\phi_0$. By the condition \eqref{eq7.5}, we obtain
\beq\label{eq7.7}
\phi_0I_0(\beta r_0)=A_2K_0(\beta r_1)=-\frac{K}{\chi}.
\eeq
By Proposition \ref{prop7.1}, we require
\beq\label{eq7.8}
\phi_0>0,\quad K<0.
\eeq
Combining this with \eqref{eq7.6} and \eqref{A2pos}, it also holds
\beq\label{eq7.9}
\phi_r(r_0)>0,\quad \phi_r(r_1)<0.
\eeq
We now consider three different cases on the region $\rho\neq 0$.

\bigskip
\noindent{\bf Case 1.} When $\frac{a\chi}{D\ve}-\beta^2=0$,
the solutions $\rho$ and $\phi$ are given as follows on $[r_0, r_1]$, 
\beq\notag
\rho(r)=\frac{c_1\chi}{\ve}\ln r+\frac{c_2\chi}{\ve}-\frac{\chi aK}{4D\ve^2}r^2+\frac{K}{\ve},
\eeq
\beq\notag
\phi(r)=c_1\ln r+c_2-\frac{aK}{4D\ve}r^2.
\eeq 
It is easy to see that 
\beq\notag
\phi_r=\frac{c_1}{r}-\frac{aK}{2D\ve}r
=\frac{1}{r}\left(-\frac{aK}{2D\ve}r^2+c_1\right).
\eeq
This contradicts \eqref{eq7.9} for any $c_1$. 
Thus, we have no nontrivial solution in this case.

\bigskip
\noindent{\bf Case 2.} When $\frac{a\chi}{D\ve}-\beta^2=-\xi^2<0$,
the solutions $\rho$ and $\phi$ are given as follows on $[r_0, r_1]$,

\beq\notag
\rho(r)=\frac{c_1\chi}{\ve}I_0(\xi r)+\frac{c_2\chi}{\ve}K_0(\xi r)+\frac{\chi aK}{D\ve^2\xi^2}+\frac{K}{\ve},
\eeq
\beq\notag
\phi(r)=c_1 I_0(\xi r)+c_2K_0(\xi r)+\frac{aK}{D\ve\xi^2}.
\eeq 

We can solve for $c_1$ and $c_2$ at $r_0$ and $r_1$ by using the fact $\phi\in C^2$  and Cramer's rule
\beq\label{eq7.15}
c_1=\frac{\tilde W_1( r_0)}{\tilde W( r_0)}=\frac{\tilde W_1( r_1)}{\tilde W( r_1)},\quad c_2=\frac{\tilde W_2( r_0)}{\tilde W( r_0)}=\frac{\tilde W_2( r_1)}{\tilde W( r_1)}
\eeq
where 
\beq\label{eq7.16}
\tilde W( r)=
\left|
\begin{array}{cc}
I_0(\xi r)&K_0(\xi r)\\
\partial_rI_0(\xi r)&\partial_r K_0(\xi r)
\end{array}
\right|
\eeq
is the Wronskian of $I_0$ and $K_0$ at $\xi r$,  and 
\beq\label{eqn7.17}
\tilde W_1(r_0)=
\left|
\begin{array}{cc}
-\frac{aK}{D\ve\xi^2}-\frac{K}{\chi}&K_0(\xi r_0)\\
\phi_0\partial_r I_0(\beta r_0 )&\partial_r K_0(\xi r_0)
\end{array}
\right|,
\quad
\tilde W_2(r_0)=
\left|
\begin{array}{cc}
I_0(\xi r_0)&-\frac{aK}{D\ve\xi^2}-\frac{K}{\chi}\\
\partial_rI_0(\xi r_0)&\phi_0\partial_r I_0(\beta r_0 )
\end{array}
\right|,
\eeq
\beq\label{eqn7.18}
\tilde W_1(r_1)=
\left|
\begin{array}{cc}
-\frac{aK}{D\ve\xi^2}-\frac{K}{\chi}&K_0(\xi r_1)\\
A_2\partial_r K_0(\beta r_1 )&\partial_r K_0(\xi r_1)
\end{array}
\right|,
\quad
\tilde W_2(r_1)=
\left|
\begin{array}{cc}
I_0(\xi r_1)&-\frac{aK}{D\ve\xi^2}-\frac{K}{\chi}\\
\partial_rI_0(\xi r_1)&A_2\partial_r K_0(\beta r_1 )
\end{array}
\right|.
\eeq
Since $I_0(r)>0$ is increasing and $K_0(r)> 0$ is decreasing for any $0< r<\infty$,  we have 
\beq\notag
\tilde W(r)<0.
\eeq
Direct computation implies 
\beq\notag
-\frac{aK}{D\ve\xi^2}-\frac{K}{\chi}=-\frac{K\beta^2}{\chi \xi^2}>0.
\eeq
Thus, $\tilde W_1(r_0)<0$ and $\tilde W_2(r_1)<0$, which implies 
\beq\label{eq7.21}
c_1>0,\quad c_2>0.
\eeq
Taking the derivative with respect to $r$, it holds
\beq\notag
\phi_r(r)=c_1\partial_rI_0(\xi r)+c_2\partial_rK_0(\xi r).
\eeq
By the fact that $\partial_rI_0$ and $\partial_rK_0$ are increasing and \eqref{eq7.21}, we find a contradiction to condition \eqref{eq7.9}, which implies that there is no nontrivial solution in this case. 

\bigskip
\noindent{\bf Case 3.} When $\frac{a\chi}{D\ve}-\beta^2=\omega^2>0$, 
the solutions $\rho$ and $\phi$ are given as follows on $[r_0, r_1]$,
\beq\notag
\rho(r)=\frac{c_1\chi}{\ve}J_0(\omega r)+\frac{c_2\chi}{\ve}Y_0(\omega r)-\frac{\chi aK}{D\ve^2\omega^2}+\frac{K}{\ve},
\eeq
\beq\notag
\phi(r)=c_1 J_0(\omega r)+c_2Y_0(\omega r)-\frac{aK}{D\ve\omega^2}.
\eeq

We can solve for $c_1$ and $c_2$ at $r_0$ and $r_1$ by using the fact $\phi\in C^2$  and Cramer's rule
\beq\label{eq7.25}
c_1=\frac{M_1( r_0)}{M( r_0)}=\frac{M_1( r_1)}{M( r_1)},\quad c_2=\frac{M_2( r_0)}{M( r_0)}=\frac{M_2( r_1)}{M( r_1)}
\eeq
where 
\beq\label{eq7.26}
 M( r)=
\left|
\begin{array}{cc}
J_0(\xi r)&Y_0(\xi r)\\
\partial_rJ_0(\xi r)&\partial_r Y_0(\xi r)
\end{array}
\right|
\eeq
is the Wronskian of $J_0$ and $Y_0$ at $\xi r$,  and 
\beq\label{eqn7.27}
M_1(r_0)=
\left|
\begin{array}{cc}
\frac{aK}{D\ve\omega^2}-\frac{K}{\chi}&Y_0(\xi r_0)\\
\phi_0\partial_r I_0(\beta r_0 )&\partial_r Y_0(\xi r_0)
\end{array}
\right|,
\quad
M_2(r_0)=
\left|
\begin{array}{cc}
J_0(\xi r_0)&\frac{aK}{D\ve\omega^2}-\frac{K}{\chi}\\
\partial_rJ_0(\xi r_0)&\phi_0\partial_r I_0(\beta r_0 )
\end{array}
\right|,
\eeq
\beq\label{eqn7.28}
M_1(r_1)=
\left|
\begin{array}{cc}
\frac{aK}{D\ve\omega^2}-\frac{K}{\chi}&Y_0(\xi r_1)\\
A_2\partial_r K_0(\beta r_1 )&\partial_r Y_0(\xi r_1)
\end{array}
\right|,
\quad
M_2(r_1)=
\left|
\begin{array}{cc}
J_0(\xi r_1)&-\frac{aK}{D\ve\xi^2}-\frac{K}{\chi}\\
\partial_rJ_0(\xi r_1)&A_2\partial_r K_0(\beta r_1 )
\end{array}
\right|.
\eeq
Direct computation implies 
\beq\notag
\frac{aK}{D\ve\omega^2}-\frac{K}{\chi}=-\frac{K\beta^2}{\chi \omega^2}>0.
\eeq
Therefore, if we require $r_0$ and $r_1$ to satisfy
\beq\label{eq7.30}
\frac{M_1( r_0)}{M( r_0)}=\frac{M_1( r_1)}{M( r_1)},\quad \frac{M_2( r_0)}{M( r_0)}=\frac{M_2( r_1)}{M( r_1)},
\eeq
and 
\beq\label{eq7.31}
c_1\partial_rJ_0(\xi r_0)+c_2\partial_rY_0(\xi r_0)>0,\quad
c_1\partial_rJ_0(\xi r_1)+c_2\partial_rY_0(\xi r_1)<0
\eeq
according the condition \eqref{eq7.9}, then there are nontrivial single bump solutions to system \eqref{steqn} in this case.

%%%%%%%%%%%%%%%%%%%%%%%%%%%%%%%%%%%
%%%%%%%%%%%%%%%%%%%%%%%%%%%%%%%%%%%

\bigskip
\begin{appendices}
\section{Energy estimates}
\label{App1}
\setcounter{equation}{0}

This section is devoted to the energy estimates of the hyperbolic-parabolic system \eqref{flow2D} with  
$$p(\rho)=\frac{\ve(\gamma-1)}{\gamma}\rho^{\gamma}.$$ 
We first define the energy of the system as follows
\beq\label{engJ}
E(t)=\int \left( \frac{1}{2}\rho |u|^2+\frac{\ve}{\gamma}\rho^{\gamma}+\frac{\chi D}{2a}|\nabla \phi|^2+\frac{\chi b}{2a}|\phi|^2-\chi \rho\phi\right)\,d{\bf x}.
\eeq
Here, ${\bf x}=(x,y)\in \R^2$, and we denote $\int\,d{\bf x}=\int_{\R^2}\,d{\bf x}$ for simplicity.

\begin{proposition}\label{propA1}
For any smooth solution $(\rho, u, \phi)$ to the system \eqref{flow2D}, we have the following energy estimate
\beq\label{enginq1}
\frac{d}{dt} E(t)=-\int\left(\frac{\chi\delta}{a}|\phi_t|^2+\alpha\rho|u|^2\right)\,d{\bf x}\leq 0.
\eeq
\end{proposition}

\pf Multiplying the second equation of system \eqref{flow2D} by $u$ and integrating over $\R^2$, we obtain
\beq\label{engpf1}
\int\left((\rho u)_t+\di (\rho u\otimes u)\right)\cdot u\,d{\bf x}
+\int \nabla p(\rho)\cdot u\,d{\bf x}
=\int \left(\chi\rho\nabla\phi\cdot u-\alpha\rho |u|^2\right)\,d{\bf x}.
\eeq
For the first term of the left side, utilizing the first equation of \eqref{flow2D}, we have
\beq\label{engpf2}
\begin{split}
\int\left((\rho u)_t+\di (\rho u\otimes u)\right)\cdot u\,d{\bf x}
=\int\left(\rho u_t+\rho u\cdot\nabla u\right)\cdot u\,d{\bf x}
=\frac{d}{dt}\int \frac{1}{2}\rho |u|^2\,d{\bf x}.
\end{split}
\eeq
For the second term of the left side, direct computation implies
\beq\notag
\begin{split}
\int \nabla p(\rho)\cdot u\,d{\bf x}
=\int\ve(\gamma-1)\rho^{\gamma-1}u\cdot \nabla \rho\,d{\bf x}
=-\int\ve(\gamma-1)\rho^{\gamma-1}(\rho_t+\rho\,\di u)\,d{\bf x}.
\end{split}
\eeq
Integrating by parts for the second term, it holds
\beq\notag
\begin{split}
\int \nabla p(\rho)\cdot u\,d{\bf x}
=-\frac{d}{dt}\int\frac{\ve(\gamma-1)}{\gamma}\rho^{\gamma}\,d{\bf x}
+\int\ve(\gamma-1)\gamma\rho^{\gamma-1}u\cdot\nabla \rho\,d{\bf x}.
\end{split}
\eeq
Hence,
\beq\label{engpf3}
\int \nabla p(\rho)\cdot u\,d{\bf x}
=\frac{d}{dt}\int\frac{\ve}{\gamma}\rho^{\gamma}\,d{\bf x}.
\eeq
For the first term of the right side, integrating by parts and using the first equation of \eqref{flow2D}, we obtain
\beq\label{engpf4}
\int\chi\rho u\cdot\nabla \phi\,d{\bf x}
=-\int\chi\di(\rho u)\phi\,d{\bf x}
=\int\chi\rho_t\phi\,d{\bf x}
=\frac{d}{dt}\int\chi\rho \phi\,d{\bf x}-\int\chi\rho \phi_t\,d{\bf x}.
\eeq
Combining \eqref{engpf2}-\eqref{engpf4} with \eqref{engpf1}, we obtain
\beq\label{engpf5}
\frac{d}{dt}\int\left(\frac{1}{2}\rho|u|^2+\frac{\ve}{\gamma}\rho^{\gamma}-\chi\rho\phi\right)\,d{\bf x}
=-\int \left(\chi\rho\phi_t+\alpha\rho |u|^2\right)\,d{\bf x}.
\eeq
Multiplying the third equation of \eqref{flow2D} by $\frac{\chi}{a}\phi_t$ and integrating by parts over $\R^2$, we obtain
\beq\label{engpf6}
\frac{d}{dt}\int\left(\frac{\chi D}{2a}|\nabla \phi|^2+\frac{\chi b}{2a}|\phi|^2\right)\,d{\bf x}
=\int\left(\chi\rho\phi_t-\frac{\chi\delta}{a}|\phi_t|^2\right)\,d{\bf x}.
\eeq
Combining \eqref{engpf5} and \eqref{engpf6}, we conclude \eqref{enginq1}, which completes the proof.
\endpf

\begin{remark}
Same estimates can be obtain on bounded region with smooth boundary and suitable boundary conditions. 
\end{remark}

\begin{proposition}\label{propA3}
Suppose that smooth solution $(\rho, u, \phi)$ to the system \eqref{flow2D} satisfies $\rho\geq 0$ and $\int\rho(x,t)\,d{\bf x}=m$ for some constant $m\geq 0$. Then the energy functional $E(t)$ is bounded from below and we have the following estimates for any $T>0$
\beq\label{enginq2}
 E^+(T)+2\int_0^T\int\left(\frac{\chi\delta}{a}|\phi_t|^2+\alpha\rho|u|^2\right)\,d{\bf x}dt\leq 2E(0)+Cm^{q_0}
\eeq
where $q_0>0$ is constant depending on $\gamma$ and
\beq\label{defE+}
E^+(t)=\int \left( \frac{1}{2}\rho |u|^2+\frac{\ve}{\gamma}\rho^{\gamma}+\frac{\chi D}{2a}|\nabla \phi|^2+\frac{\chi b}{2a}|\phi|^2\right)\,d{\bf x}.
\eeq
\end{proposition}

\pf
By Proposition \ref{propA1} and the relation $E^+(t)=E(t)+\int\chi\rho\phi\,d{\bf x}$, we only need to estimate the last term. We divide our proof into two cases.

\medskip
\noindent{\bf Case 1.} When $\gamma> 2$, we obtain by the H$\ddot{\mbox{o}}$lder inequality and the interpolation inequality 
\beq\notag
\begin{split}
\left|\int\chi\rho\phi\,d{\bf x}\right|
\leq C\|\rho\|_{2}\|\phi\|_{2}
\leq C\|\rho\|_{1}^{\eta}\|\rho\|_{\gamma}^{1-\eta}\|\phi\|_{2}
\leq Cm^{\eta}\big(E^+(t)\big)^{\frac{1-\eta}{\gamma}+\frac12}
\end{split}
\eeq
where $0<\eta<1$ and  
$$
 \eta=\frac{\gamma-2}{2(\gamma-1)}.
$$
Direct computation implies 
\beq\notag
\frac{1-\eta}{\gamma}+\frac12=\frac{1}{2(\gamma-1)}+\frac12<1.
\eeq
Utilizing the Young inequality, it holds
\beq\label{engpf71}
\begin{split}
\left|\int\chi\rho\phi\,d{\bf x}\right|
\leq Cm^{q_0}+\frac{1}{2} E^{+}(t),
\end{split}
\eeq
where 
$$
q_0=\frac{2\eta\gamma}{\gamma-1+\eta}.
$$
It is not hard to see from this inequality that $E(t)$ is bounded from below
\beq\notag
E(t)=E^{+}(t)-\int\chi\rho\phi\,d{\bf x}\geq \frac12E^{+}(t)-Cm^{q_0}.
\eeq
Integrating \eqref{enginq1} over $(0,T)$ and using \eqref{engpf71}, we conclude the estimate \eqref{enginq2} for this case. 

\medskip
\noindent{\bf Case 2.} When $\gamma=2$, we obtain by the H$\ddot{\mbox{o}}$lder inequality, the interpolation inequality and the Gagliardo–Nirenberg inequality
\beq\notag
\begin{split}
\left|\int\chi\rho\phi\,d{\bf x}\right|
\leq C\|\rho\|_{l}\|\phi\|_{l'}
\leq C\|\rho\|_{1}^{\eta}\|\rho\|_{2}^{1-\eta}\|\phi\|_{2}^{\zeta}\|\nabla \phi\|_{2}^{1-\zeta}
\leq Cm^{\eta}\big(E^+(t)\big)^{\frac{2-\eta}{2}}
\end{split}
\eeq
where $1<l<2$, $l'>2$, $0<\eta,\zeta<1$ and  
$$
\frac{1}{l}+\frac{1}{l'}=1,\quad \eta=\frac{2-l}{l},\quad \zeta=\frac{2}{l'}.
$$
It is easy to see 
\beq\notag
\frac{2-\eta}{2}<1.
\eeq
Utilizing the Young inequality, it holds
\beq\label{engpf72}
\begin{split}
\left|\int\chi\rho\phi\,d{\bf x}\right|
\leq Cm^{2}+\frac{1}{2} E^{+}(t).
\end{split}
\eeq
It is not hard to see from this inequality that $E(t)$ is bounded from below
\beq\notag
E(t)=E^{+}(t)-\int\chi\rho\phi\,d{\bf x}\geq \frac12E^{+}(t)-Cm^{2}.
\eeq
Integrating \eqref{enginq1} over $(0,T)$ and using \eqref{engpf72}, we conclude the estimate \eqref{enginq2} for this case with $q_0=2$.

\medskip
\noindent{\bf Case 3.} When $1<\gamma<2$, we obtain by the H$\ddot{\mbox{o}}$lder inequality, the interpolation inequality and the Gagliardo–Nirenberg inequality
\beq\notag
\begin{split}
\left|\int\chi\rho\phi\,d{\bf x}\right|
\leq C\|\rho\|_{l}\|\phi\|_{l'}
\leq C\|\rho\|_{1}^{\eta}\|\rho\|_{\gamma}^{1-\eta}\|\phi\|_{2}^{\zeta}\|\nabla \phi\|_{2}^{1-\zeta}
\leq Cm^{\eta}\big(E^+(t)\big)^{\frac{1-\eta}{\gamma}+\frac12}
\end{split}
\eeq
where $1<l<\frac{2}{3-\gamma}<\gamma<2$, $l'>2$, $0<\eta,\zeta<1$ and  
$$
\frac{1}{l}+\frac{1}{l'}=1,\quad \eta=\frac{\gamma-l}{l(\gamma-1)},\quad \zeta=\frac{2}{l'}.
$$
Direct computation implies 
\beq\notag
\frac{1-\eta}{\gamma}+\frac12=\frac{l-1}{l(\gamma-1)}+\frac12
<\frac{\frac{2}{3-\gamma}-1}{\frac{2}{3-\gamma}(\gamma-1)}+\frac12=1.
\eeq
Utilizing the Young inequality, it holds
\beq\label{engpf7}
\begin{split}
\left|\int\chi\rho\phi\,d{\bf x}\right|
\leq Cm^{q_0}+\frac{1}{2} E^{+}(t)
\end{split}
\eeq
where 
$$
q_0=\frac{2\eta\gamma}{\gamma-1+\eta}.
$$
It is not hard to see from this inequality that $E(t)$ is bounded from below
\beq\notag
E(t)=E^{+}(t)-\int\chi\rho\phi\,d{\bf x}\geq \frac12E^{+}(t)-Cm^{q_0}.
\eeq
Integrating \eqref{enginq1} over $(0,T)$ and using \eqref{engpf7}, we conclude the estimate \eqref{enginq2}, which completes the proof. 

\endpf

\end{appendices}

%%%%%%%%%%%%%%%%%%%%%%%%%%%%%%%%%%%
%%%%%%%%%%%%%%%%%%%%%%%%%%%%%%%%%%%

%\bibliographystyle{plain}
%\bibliography{bib.bib}

\begin{thebibliography}{99}

\bibitem{2016Paper}
F. Berthelin, D. Chiron and M. Ribot, 
{\em Stationary solutions with vacuum for a one-dimensional chemotaxis
  model with non-linear pressure}, Comm. Math. Sci., 14(1), 147--186, 2016.
  
 \bibitem{CCWWZ20} 
J. Carrillo, X. Chen, Q. Wang, Z. Wang and L. Zhang, {\em Phase transitions and bump solutions of the Keller-Segel model with volume exclusion}, SIAM J. Appl. Math., 80, 232--261, 2020. 
  
  \bibitem{CHS23}
  T. Crin-Barat, Q.Y. He and L.Y. Shou, {\em The hyperbolic-parabolic chemotaxis system modelling vasculogenesis: global dynamics and relaxation limit toward a Keller-Segel model}, SIAM J. Math. Anal., 55(5), 4445--4492, 2023.
  
  \bibitem{Russo12}
C. Di Russo, {\em Analysis and numerical approximations of hydrodynamical models of biological movements}, Rend. Mat. Appl., 32(3-4), 117--367, 2012.

\bibitem{Russo13}
C. Di Russo and A. Sepe, {\em Existence and asymptotic behavior of solutions to a quasi-linear
hyperbolic-parabolic model of vasculogenesis}, SIAM J. Math. Anal., 45(2), 748--776, 2013.

\bibitem{chemotaxis}
M. Eisenbach, A. Tamada, G.M. Omann, J.E. Segall, R.A. Firtel, R. Meili,
  D. Gutnick, M. Varon, J.W. Lengeler and F. Murakami,
 {\em Chemotaxis}, World Scientific Publishing Company, 2004.

\bibitem{Gamba03}
A. Gamba, D. Ambrosi, A. Coniglio, A de Candia, S. Di Talia, E. Giraudo, et al., {\em Percolation, morphogenesis, and Burgers dynamics in blood vessels formation}, Phys. Rev. Lett., 90(11), 118101, 2003,

\bibitem{PDEmodels}
T. Hillen and K. Painter,
{\em A user’s guide to PDE models for chemotaxis}, J. Math. Bio., 58, 183--217, 2008.


\bibitem{HPWZ21} 
G.Y. Hong, H.Y. Peng, Z.A. Wang and C.J. Zhu, {\em Nonlinear stability of phase transition steady states to a hyperbolic–parabolic system modeling vascular networks}, J. London Math. Soc., 103, 1480--1514, 2021. 


\bibitem{LW24}
M.Q. Liu and Z.G. Wu, {\em Large time behavior of a hyperbolic-parabolic model of vasculogenesis},
Discrete Contin. Dyn. Sys. Series B, 29(2), 777--795, 2024. 


\bibitem{LPW221}
 Q.Q. Liu, H.Y. Peng and Z.A. Wang, {\em Asymptotic stability of diffusion waves of a quasi-linear hyperbolic-parabolic model for vasculogenesis}, SIAM J. Math. Anal., 54, 1313--1346, 2022. 
 
 \bibitem{LPW222}
Q.Q. Liu, H.Y. Peng and Z.A. Wang, {\em Convergence to nonlinear diffusion waves for a hyperbolic- parabolic chemotaxis system modelling vasculogenesis}, J. Differ. Equ., 314, 251--286, 2022. 

\bibitem{PZ23}
H.Y. Peng and K. Zhao, {\em On a hyperbolic-parabolic chemotaxis system}, Math. Biosci. Engineering, 20(5), 7802--7827, 2023.


\bibitem{Serini03}
G. Serini, D. Ambrosi, E. Giraudo, A. Gamba, L. Preziosi and F. Bussolino, {\em Modeling the early stages of vascular network assembly}, The EMBO Journal, 22(8), 1771--1779, 2003.

\bibitem{Roman00}
B.L. Roman and B.M. Weinstein, {\em Building the vertebrate vasculature: research is going swimmingly}, Bioessays, 22(10), 882--893, 2000.


\end{thebibliography}

\end{document}